\documentclass[10pt,epsfig]{article}
\usepackage{amsbsy,amssymb,amscd,amsfonts,latexsym,amstext,delarray,
amsmath,epsfig}  \headsep=15pt

\input xypic

\newtheorem{thm}{Theorem}[section]
\newtheorem{prop}[thm]{Proposition}

\newtheorem{lem}[thm]{Lemma}

\newtheorem{defn}[thm]{Definition}

\numberwithin{equation}{section}

\def\R{{\mathbb R}}

\def\H{{\mathbb H}}
\def\m{{\mathfrak m}}

\def\O{{\mathcal O}}
\def\mX{{\mathfrak X}}
\def\mK{{\mathbb K}}
\def\Sp{{\rm{Spec}}}

\def\Q{{\mathbb Q}}
\def\C{{\mathbb C}}
\def\Z{{\mathbb Z}}
\def\N{{\mathbb N}}

\def\P{{\mathbb P}}

\def\PSL{{\rm PSL}}
\def\PGL{{\rm PGL}}

\def\Tr{{\rm Tr}}

\def\cW{{\mathcal W}}
\def\cP{{\mathcal P}}

\def\cS{{\mathcal S}}

\def\cH{{\mathcal H}}

\def\cV{{\mathcal V}}

\newcommand{\ie}{{\it i.e.\/}\ }

\newcommand{\cf}{{\it cf.\/}\ }

\newcommand{\resp}{{\it resp.\/}\ }

\begin{document}
\title{New perspectives in Arakelov geometry}

\author{Caterina Consani and Matilde Marcolli}
\date{}
\maketitle

\bigskip

\begin{abstract}

In this paper we give a unified description of the archimedean
and the totally split degenerate fibers of an arithmetic surface,
using operator algebras and Connes' theory of spectral triples in
noncommutative geometry.

\end{abstract}

\section{Introduction.}

The aim of this article is to report on some of the recent results
obtained by the two authors on a non-commutative interpretation of
the totally split degenerate fibers of an arithmetic surface. The
results stated in the first part of the paper were announced in a
talk given by the first author at the VII meeting of the CNTA in
May 2002 at Montreal. The material presented in this note is based
on the papers \cite{CM} and \cite{CM1}.
\smallskip

Let $\mathcal X$ be an arithmetic surface defined over $\Sp(\Z)$
(or over $\Sp({\mathcal O}_{\mK})$, for a number field $\mK$),
having the smooth algebraic curve $X_{/\Q}$ as its generic fiber.
It is well known that, as a Riemann surface, $X_{/\C}$ admits
always a uniformization by means of a Schottky group $\Gamma$. In
fact, the presence of this uniformization plays a fundamental role
in the theory of the ``fiber at infinity'' of $\mathcal X$
described in \cite{Man} and \cite{CM}. In analogy to Mumford's
p-adic uniformization of algebraic curves (\cf\cite{Mum}), the
Riemann surface $X_{/\C}$ can be interpreted as the boundary at
infinity of a 3-manifold $\mX_\Gamma$ defined as the quotient of
the real hyperbolic 3-space $\H^3$ by the action of the Schottky
group $\Gamma$. The space $\mX_\Gamma$ contains in its interior an
infinite link of bounded geodesics. Manin gave an expression for
the Arakelov Green function on $X_{/\C}$ in terms of
configurations of geodesics in $\mX_\Gamma$, thus interpreting
this tangle as the dual graph $\mathcal G$ of the ``closed fiber
at infinity'' of $\mathcal X$.

In the first part of this paper we concentrate on Manin's
description of such dual graph and we exhibit the suspension flow
$\mathcal S_T$ of a dynamical system $T$ on the limit set of the
Schottky group $\Gamma$ as a geometric model for the dual graph
${\mathcal G}$. In particular, the first cohomology group of $\mathcal
S_T$ determines a model of the first cohomology of the dual graph of
the ``fiber at infinity''. Furthermore, the first (co)homology
group of $\mathcal S_T$ carries a natural filtration.

A crucial feature of this construction is the fact (proved in \S 5
of \cite{CM}) that this dynamical cohomology contains a subspace
isomorphic to the Archimedean cohomology of \cite{KC}: the group
of invariants for the action of the local monodromy ``at
infinity''. This space has also the correct geometric properties,
in order to be interpreted in terms of a cohomology theory
associated to a degeneration in a ``neighborhood of arithmetic
infinity''. Our result identifies such space with a distinguished
subspace of the cohomology of a topological space constructed in
terms of geodesics in $\mX_\Gamma$, our geometrc model of the dual
graph. Under this identification, the graded structure associated
to the filtration on the (co)homology of $\mathcal S_T$
corresponds to the graded structure given by Tate twists on the
Archimedean cohomology of \cite{KC}.

The Cuntz-Krieger algebra ${\mathcal O}_A$ associated to the shift $T$
describes the ``ring of functions'' on a noncommutative space, which
represents the quotient of the limit set $\Lambda_\Gamma$ of the
Schottky group, by the action of $\Gamma$. In terms of the geometry of
the fiber at
arithmetic infinity, this space can be thought of as the set of
components of the special fiber, or equivalently the vertices of the
dual graph ${\mathcal G}$, whereas the quotient
$\Lambda_\Gamma\times_\Gamma \Lambda_\Gamma$ gives the edges of
${\mathcal G}$.
The algebra ${\mathcal O}_A$ carries a refined information
on the action of the Schottky group $\Gamma$ on its limit set.
In particular, we construct a {\em spectral triple} for this algebra.

In non-commutative geometry, the notion of a spectral triple
provides the correct generalization of the classical structure of
a Riemannian manifold. The two notions agree on a commutative
space. In the usual context of Riemannian geometry, the definition
of the infinitesimal element $ds$ on a smooth spin manifold can be
expressed in terms of the inverse of the classical Dirac operator
$D$. This is the key remark that motivates the theory of spectral
triples. In particular, the geodesic distance between two points
on the manifold is defined in terms of $D^{-1}$ (\cf \cite{Co94}
\S VI). The spectral triple that describes a classical Riemannian
spin manifold is $(A,H,D)$, where $A$ is the algebra of complex
valued smooth functions on the manifold, $H$ is the Hilbert space
of square integrable spinor sections, and $D$ is the classical
Dirac operator (a square root of the Laplacian). These data
determine completely and uniquely the Riemannian geometry on the
manifold. It turns out that, when expressed in this form, the
notion of spectral triple extends to more general non-commutative
spaces, where the data $(A,H,D)$ consist of a ${\rm C}^*$-algebra
$A$ (or more generally of a smooth subalgebra of a ${\rm
C}^*$-algebra) with a representation as bounded operators on a
Hilbert space $H$, and an operator $D$ on $H$ that verifies the
main properties of a Dirac operator. The notion of smoothness
is determined by $D$: the smooth elements of $A$ are defined
by the intersection of domains of powers of the derivation given
by commutator with $|D|$.
The basic geometric structure encoded by the theory of spectral
triples is Riemannian geometry, but in more refined cases, such as
K\"ahler geometry, the additional structure can be easily encoded as
additional symmetries.

In these constructions, the Dirac operator $D$ is obtained from
the grading operator associated to a filtration on the cochains of
the complex that computes the dynamical cohomology. The induced
operator on the subspace identified with the Archimedean
cohomology agrees with the ``logarithm of Frobenius'' of \cite{KC}
and \cite{Den}.

This structure further enriches the geometric interpretation of
the Archimedean cohomology, giving it the meaning of spinors on a
noncommutative manifold, with the logarithm of Frobenius
introduced in \cite{Den} in the role of the Dirac operator.

An advantage of this construction is that a completely analogous
formulation exists in the case of Mumford curves. This provides a
unified description of the archimedean and totally split
degenerate fibers of an arithmetic surface.

Let $p$ be a finite prime where ${\mathcal X}$ has totally split
degenerate reduction. Then, the completion $\hat X_p$ at $p$ of
the generic fiber of $\mX$ is a split-degenerate stable curve over
$\Q_p$ (also called a Mumford curve) uniformized by the action of
a $p$-adic Schottky group $\Gamma$. The dual graph of the
reduction of $\hat X_p$ coincides with a finite graph obtained as
the quotient of a tree $\Delta_\Gamma$ by the action of $\Gamma$.

The curve $\hat X_p$ is holomorphically isomorphic to a quotient
of a subset of the ends of the Bruhat-Tits tree $\Delta$ of $\Q_p$
by the action of $\Gamma$. Thus, in this setting, the Bruhat-Tits
tree at $p$ replaces the hyperbolic space $\H^3$ ``at infinity'',
and the analog of the tangle of bounded geodesics in $\mX_\Gamma$
is played by doubly infinite walks in $\Delta_\Gamma/\Gamma$.

In analogy with the archimedean construction, we define the system
$(\mathcal W(\Delta_\Gamma/\Gamma),T)$ where $T$ is an invertible
shift map on the set $\mathcal W(\Delta_\Gamma/\Gamma)$ of
doubly-infinite walks on the graph $\Delta_\Gamma/\Gamma$. The
first cohomology group $H^1(\mathcal
W(\Delta_\Gamma/\Gamma)_T,\Z)$ of the mapping torus $\mathcal
W(\Delta_\Gamma/\Gamma)_T$ of $T$ inherits a natural filtration
using which we introduce a dynamical cohomology group. We again
have a Cuntz-Krieger graph algebra  ${\rm C}^*(\Delta_\Gamma/\Gamma)$
and we can construct a spectral triple as in the case at infinity,
where again the Dirac operator is related to the grading operator
$\Phi$ that computes the local factor as a regularized determinant, as
in \cite{Den2}, \cite{Den3}.
In \cite{CM1}, we also suggested a possible way of extending such
construction to places that are not of split degenerate reduction,
inspired by the ``foam space'' construction of \cite{ManC} and
\cite{CMZ}.

\medskip

\noindent{\bf Acknowledgment.} This paper was partly written
during visits of the second author to Florida State University and
University of Toronto. We thank these institutions for the
hospitality. The first author is partially supported by NSERC
grant 72016789. The second author is partially supported by the
Humboldt Foundation and the German Government (Sofja Kovalevskaya
Award). We thank Alain Connes for many extremely helpful 
discussions. 

\section{Notation.}

Throughout this paper we will denote by $K$ one among the
following fields: (a) the complex numbers $\C$, (b) a finite
extension of $\Q_p$. When (b) occurs, we write $\O_K$ for the ring
of integers of $K$, $\m\subset \O_K$ for the maximal ideal and
$\pi\in\m$ for a uniformizer (\ie $\m = (\pi)$). We also denote by
$k$ the residue classes field $k = \O/\m$.\smallskip

We denote by $\H^3$ the three-dimensional real hyperbolic space
\ie the quotient $$\H^3 ={\rm SU}(2)\backslash {\rm PGL}(2,\C).$$
This space can also be described as the upper half space $\H^3
\simeq \C\times\R^+$ endowed with the hyperbolic metric. The group
$\PSL(2,\C)$ acts on $\H^3$ by isometries. The complex projective
line $\P^1(\C)$ is identified with the conformal boundary at
infinity of $\H^3$ and the action of $\PSL(2,\C)$ on $\H^3$
extends to an action on $\overline \H^3:=\H^3\cup \P^1(\C)$. The
group $\PSL(2,\C)$ acts on $\P^1(\C)$ by fractional linear
transformations.\smallskip

For an integer $g \ge 1$, a Schottky group of rank $g$ is a
discrete subgroup $\Gamma\subset \PSL (2,K)$, which is purely
loxodromic and isomorphic to a free group of rank $g$.
We denote by $\Lambda_\Gamma$ the limit set of the action of
$\Gamma$. One sees that $\Lambda_\Gamma$ is contained in
$\P^1(K)$. This set can also be described as the closure of the
set of the attractive and repelling fixed points $z^{\pm}(g)$ of
the loxodromic elements $g\in \Gamma$. In the case $g=1$ the limit
set consists of two points, but for $g\geq 2$ the limit set is
usually a fractal of some Hausdorff dimension $0\leq \delta
=\dim_H (\Lambda_\Gamma) < 2$.
We denote by $\Omega_\Gamma = \Omega_\Gamma(K)$ the domain of
discontinuity of $\Gamma$, that is, the complement of
$\Lambda_\Gamma$ in $\P^1(K)$.

When $K = \C$, the quotient space ${\mathfrak X}_\Gamma := \H^3/\Gamma$ is
topologically a handlebody of genus $g$, and the quotient $X_{/\C} =
\Omega_\Gamma /\Gamma$ is a
Riemann surface of genus $g$. The covering $\Omega_\Gamma \to
X_{/\C}$ is called a Schottky uniformization of $X_{/\C}$. Every
complex Riemann surface $X_{/\C}$ admits a Schottky
uniformization. The handlebody ${\mathfrak X}_\Gamma$ can be
compactified by adding the conformal boundary at infinity
$X_{/\C}$ to obtain $\overline{\mathfrak X}_\Gamma := {\mathfrak
X}_\Gamma\cup X_{/\C}= (\H^3 \cup
\Omega_\Gamma)/\Gamma$.\smallskip

A directed graph $E$ consists of data $E=(E^0, E^1, E^1_+, r, s,
\iota)$, where $E^0$ is the set of vertices, $E^1$ is the set of
oriented edges $w=\{ e, \epsilon \}$, where $e$ is an edge of the
graph and $\epsilon =\pm 1$ is a choice of orientation. The set
$E^1_+$ consists of a choice of orientation for each edge, namely
one element in each pair $\{ e, \pm 1 \}$. The maps $r,s: E^1 \to
E^0$ are the range and source maps, and $\iota$ is the involution
on $E^1$ defined by $\iota(w)=\{ e, -\epsilon \}$.

A directed graph is finite if $E^0$ and $E^1$ are finite sets. It
is locally finite if each vertex emits and receives at most
finitely many oriented edges in $E^1$. A vertex $v$ in a directed
graph is a sink if there is no edge in $E^1_+$ with source $v$. A
juxtaposition of oriented edges $w_1 w_2$ is said to be admissible
if $w_2\neq \iota(w_1)$ and $r(w_1)=s(w_2)$. A (finite, infinite,
doubly infinite) walk in a directed graph $E$ is an admissible
(finite, infinite, doubly infinite) sequence of elements in $E^1$.

We denote by $\cW^n(E)$ the set of walks of length $n$, by
$\cW^*(E)=\cup_n \cW^n(E)$, by $\cW^+(E)$ the set of infinite
walks, and by $\cW(E)$ the set of doubly infinite walks. A
directed graph is a directed tree if, for any two vertices, there
exists a unique walk in $\cW^*(E)$ connecting them.

The edge matrix $A_+$ of a locally finite (or row finite) directed
graph is an $(\# E^1_+)\times (\# E^1_+)$ (possibly infinite)
matrix. The entries $A_+(w_i,w_j)$ satisfy $A_+(w_i,w_j)=1$ if
$w_i w_j$ is an admissible path, and $A_+(w_i,w_j)=0$ otherwise.
The directed edge matrix of $E$ is a $\# E^1 \times \# E^1$
(possibly infinite) matrix with entries $A(w_i,w_j)=1$ if $w_i
w_j$ is an admissible walk and $A(w_i,w_j)=0$ otherwise.\smallskip

Even when not explicitly stated, all Hilbert spaces and algebras
of operators we consider will be separable, \ie they admit a dense
(in the norm topology) countable subset.

\section{A dynamical theory at infinity.}\label{6}

Let $\Gamma\subset \PSL (2,\C)$ be a Schottky group of rank $g \ge
2$. Given a choice of a set of generators $\{ g_i \}_{i=1}^g$ for
$\Gamma$, there is a bijection between the elements of $\Gamma$
and the set of all admissible walks in the Cayley graph of
$\Gamma$, namely reduced words in the $\{ g_i \}_{i=1}^{2g}$,
where we use the notation $g_{i+g} := g_{i}^{-1}$, for
$i=1,\ldots, g$.

In the following we consider the sets ${\mathcal S^+}$ and
${\mathcal S}$ of \resp right-infinite, doubly infinite admissible
sequences in the $\{ g_i \}_{i=1}^{2g}$:
\begin{equation}\label{S+} {\mathcal S}^+= \{ a_0 a_1 \ldots
a_\ell \ldots  \, \, | a_i \in \{ g_i \}_{i=1}^{2g}, \, \, a_{i+1}
\neq a_i^{-1}, \forall i\in \N \},
\end{equation}
\begin{equation}\label{S} {\mathcal S}= \{\ldots a_{-m} \ldots a_{-1}
a_0 a_1 \ldots a_\ell \ldots  \, \, | a_i \in \{ g_i
\}_{i=1}^{2g}, \, \, a_{i+1} \neq a_i^{-1}, \forall i\in \Z \}.
\end{equation}
The admissibility condition simply means that we only allow
``reduced'' words in the generators, without cancellations.

On the space ${\mathcal S}$ we consider the topology generated by
the sets $W^s(x,\ell) =\{ y\in {\mathcal S} | x_k = y_k, k\geq
\ell \}$, and $W^u(x,\ell)=\{ y\in {\mathcal S} | x_k = y_k, k\leq
\ell \}$ for $x\in {\mathcal S}$ and $\ell \in \Z$.

There is a two-sided shift operator $T$ acting on ${\mathcal S}$
as the map
\begin{equation}\label{shift}
\begin{array}{rcccccccccl} T(& \ldots & a_{-m} & \ldots & a_{-1} & a_0
& a_1 & \ldots & a_{\ell} & \ldots &) = \\ & \ldots & a_{-m+1} &
\ldots & a_{0} & a_1 & a_2 & \ldots & a_{\ell+1} & \ldots &
\end{array} \end{equation}

\begin{defn}
A subshift of finite type $(\cS_A,T)$ consists of all doubly
infinite sequences in the elements of a given finite set $W$
(alphabet) with the admissibility condition specified by a $\# W
\times \# W$ elementary matrix,
$$ \cS_A =\{ (w_k)_{k\in \Z}: w_k \in W, A(w_k,w_{k+1})=1 \}, $$
and with the action of the invertible shift $(Tw)_k =w_{k+1}$.
\end{defn}

\begin{lem} \label{S-Lamb2}
The space ${\mathcal S}$ can be identified with the subshift of
finite type ${\mathcal S}_A$ with the symmetric $2g\times 2g$
matrix $A$ given by the directed edge matrix of the Cayley graph
of $\Gamma$.
\end{lem}

The two-sided shift operator $T$ on ${\mathcal S}$ of
\eqref{shift} decomposes ${\mathcal S}$ in a product of expanding
and contracting directions, so that $( {\mathcal S}, T)$ is a
Smale space.

The following topological space is defined in terms of the Smale
space $({\mathcal S},T)$ and will be considered as a geometric
realization of the  ``dual graph'' associated to the fiber at
arithmetic infinity of the arithmetic surface ${\mathcal X}$.

\begin{defn} The mapping torus (suspension flow) of the dynamical
system $( {\mathcal S}, T)$ is defined as
\begin{equation} \label{suspensionT}  {\mathcal S}_T := {\mathcal S}
\times [0,1] / (x,0)\sim (Tx,1) \end{equation}
\end{defn}\medskip

The space ${\mathcal S}_T$ is a very natural space associated to the 
noncommutative space
\begin{equation}\label{quotientNC}
\Lambda_\Gamma \times_\Gamma \Lambda_\Gamma \simeq {\mathcal S}/\Z,
\end{equation}
with $\Z$ acting via the invertible shift $T$ of \eqref{shift}, namely
the {\em homotopy quotient} (\cf \cite{BC} \cite{Co-tr}),
\begin{equation}\label{htpyquot}
{\mathcal S}_T = {\mathcal S}\times_\Z \R.
\end{equation}
Namely, it is a commutative space that provides, up to homotopy, a
geometric model for \eqref{quotientNC}, where
the noncommutative space \eqref{quotientNC} can be identified with the
quotient space of a foliation \eqref{htpyquot} whose generic leaf is
contractible (a copy of $\R$).

\subsection{(Co)homology of ${\mathcal S}_T$.}

In this paragraph we give an explicit description of the
(co)homology $H^1({\mathcal S}_T,\Z)$.\smallskip

The shift $T$ acting on ${\mathcal S}$ induces an automorphism of
the $C^*$--algebra of continuous functions $C({\mathcal S})$. With
an abuse of notation we still denote it by $T$. Consider the
crossed product $C^*$--algebra $C({\mathcal S}) \rtimes_T \Z$.
This is a suitable norm completion of $C({\mathcal S}) [T,
T^{-1}]$ with product $(V*W)_k = \sum_{r\in\Z} V_k \cdot (T^r
W_{r+k})$, for $V=\sum_k V_k T^k$, $W=\sum_k W_k T^k$, and
$V*W=\sum_k (V*W)_k T^k$.
The $K$--theory group $K_0({\rm C}({\mathcal S}) \rtimes_T \Z)$ is
described by the co--invariants of the action of $T$ (\cf
\cite{BoHa}, \cite{PaTu}).

\begin{thm}
The cohomology $H^1({\mathcal S}_T)$ satisfies the following
properties
\begin{enumerate}
\item There is an identification of $H^1({\mathcal S}_T,\Z)$ with
the $K_0$-group of the crossed product ${\rm C}^*$-algebra for the
action of $T$ on ${\mathcal S}$,
\begin{equation}  H^1({\mathcal S}_T,\Z) \cong K_0({\rm C}({\mathcal
S}) \rtimes_T \Z).  \label{H1K0} \end{equation} \item The
identification \eqref{H1K0} endows $H^1({\mathcal S}_T,\Z)$ with a
filtration by free abelian groups $F_0 \hookrightarrow F_1
\hookrightarrow \cdots F_n\hookrightarrow \cdots$, with ${\rm
rank} F_0=2g$ and ${\rm rank} F_n = 2g(2g-1)^{n-1}(2g-2) +1$, for
$n\geq 1$, so that
$$ H^1({\mathcal S}_T,\Z)= \varinjlim_{n} F_n. $$
\item This filtration is induced by a filtration ${\mathcal P}_n$ on
the locally constant functions ${\rm C}({\mathcal S}^+,\Z)$ which
depend on ``future coordinates'', with ${\mathcal P}_n$ given by
locally constant functions that depend only on the first $n+1$
coordinates. The Pimsner--Voiculescu exact sequence for the $K$-theory
\eqref{H1K0},
\begin{equation}\label{PV}
0 \to \Z \to {\rm C}({\mathcal S},\Z) \stackrel{\delta=1-T}{\longrightarrow}
{\rm C}({\mathcal S},\Z) \to H^1({\mathcal S}_T,\Z) \to 0.
\end{equation}
\item The cohomology $H^1({\mathcal S}_T,\kappa)=H^1({\mathcal
S}_T,\Z)\otimes \kappa$, for $\kappa=\R$ or $\C$, is computed by
\begin{equation}\label{PVC}
0 \to \C \to {\mathcal P}\otimes \kappa  \stackrel{\delta}{\longrightarrow}
{\mathcal P}\otimes \kappa \to H^1({\mathcal S}_T,\C) \to 0.
\end{equation}
The vector space ${\mathcal P}\otimes \kappa$, for
${\mathcal P}={\rm C}({\mathcal S}^+,\Z)=\varinjlim_n {\mathcal P}_n$,
admits a Hilbert space completion ${\mathcal L}= L^2(\Lambda_\Gamma,
\mu)$, where $\mu$ is the Patterson--Sullivan measure on the limit set
$\Lambda_\Gamma$ (\cf \cite{Sull}) satisfying
\begin{equation}\label{PSmeas}
(\gamma^* d\mu)(x)= |\gamma^\prime (x)|^{\delta_H} \, d\mu(x), \ \ \
\forall \gamma\in \Gamma.
\end{equation}
\end{enumerate}
\label{thm-homologydyn}
\end{thm}

Notice that the adjoint $\delta^*$ in the $L^2$-inner product of the
coboundary $\delta$ of the
cohomology $H^1({\mathcal S}_T)$ is an important operator
associated to the dynamics of the (one sided) shift $T$ on the limit
set $\Lambda_\Gamma$, namely, the operator $1-{\mathcal R}$, where
${\mathcal R}$ is the Perron--Frobenius operator of $T$.
This is the analog of the Gauss--Kuzmin operator studied in
\cite{ManMar} \cite{Mar} in the case of modular curves.

\smallskip

For simplicity of notation, in the following we will use the same
notation ${\mathcal P}$
and ${\mathcal P}_n$ for the $\Z$-modules and for the $\kappa$-vector
spaces. It will be clear from the context which one we refer to.
We identify the vector spaces ${\mathcal P}_n$ with finite dimensional
subspaces of $L^2(\Lambda_\Gamma,\mu)$, by identifying locally
constant functions on ${\mathcal S}^+$ with locally constant functions
on $\Lambda_\Gamma$.

The following result computes the first homology of ${\mathcal
S}_T$.

\begin{prop}
The homology group $H_1({\mathcal S}_T,\Z)$ has a filtration by
free abelian groups ${\mathcal K}_N$,
\begin{equation}\label{hom-lim-dir}
 H_1({\mathcal S}_T,\Z) = \varinjlim_N {\mathcal K}_N,
\end{equation}
with
$$ K_N = {\rm rank}( {\mathcal K}_N ) = \left\{ \begin{array}{lr}
(2g-1)^N +1 & N \text{ even } \\ (2g-1)^N + (2g-1) & N \text{ odd
}
\end{array}\right. $$
The group $H_1({\mathcal S}_T,\Z)$ can also be written as
$$ H_1({\mathcal S}_T,\Z) = \oplus_{N=0}^\infty {\mathcal R}_N $$
where ${\mathcal R}_n$ is a free abelian group of ranks $R_1=2g$
and
$$
 R_N = {\rm rank}({\mathcal R}_N) =\frac{1}{N} \sum_{d|N} \mu(d)\,
(2g-1)^{N/d},
$$
for $N>1$, with $\mu$ the M\"obius function. This is isomorphic to
a free abelian group on countably many generators.
\label{homology-dyn}
\end{prop}

Thus, the $\Z$-modules $F_n$ are obtained as quotients $\cP_n/\delta
\cP_{n-1}$, for $\delta(f)=f-f\circ T$, where $\cP=\cup_n \cP_n$
is the module of continuous $\Z$-valued functions on ${\mathcal
S}$ that depend only on ``future coordinates''. The $\Z$-modules
${\mathcal K}_N$ are generated by all admissible words $a_0\ldots
a_N$ such that the word $a_N a_0$ is also admissible.
Combining Theorem \ref{thm-homologydyn} with Proposition
\ref{homology-dyn}, we can compute explicitly the pairing of
homology and cohomology for ${\mathcal S}_T$.

\begin{prop}
Let $F_n$ and ${\mathcal K}_N$ be the filtrations defined,
respectively, in Theorem \ref{thm-homologydyn} and Proposition
\ref{homology-dyn}. There is a pairing
\begin{equation}\label{pairing1}
\langle \cdot, \cdot \rangle : F_n \times {\mathcal K}_N \to \Z
\ \ \ \ \  \langle [f], x \rangle = N \cdot f(\bar x),
\end{equation} with $x=a_0 \ldots a_N$. Here the representative
$f\in [f]$ is a function that depends on the first $n+1$ terms
$a_0\ldots a_n$ of sequences in ${\mathcal S}$, and $\bar x$ is
the truncation of the periodic sequence $\overline{a_0\ldots a_N}$
after the first $n$ terms. This pairing descends to the direct
limits of the filtrations, where it agrees with the classical
cohomology/homology pairing
\begin{equation}\label{pairing2}
\langle \cdot, \cdot \rangle : H^1({\mathcal S}_T,\Z) \times
H_1({\mathcal S}_T,\Z) \to \Z.
\end{equation}
\label{H1pairing}
\end{prop}

\subsection{Dynamical (co)homology.}

We define the dynamical cohomology $H^1_{dyn}$ as the graded
vector space given by the sum of the graded pieces of the
filtration of $H^1({\mathcal S}_T)$, introduced in Theorem
\ref{thm-homologydyn}. These graded pieces ${\rm Gr}_n$ are
considered with coefficients in the $n$-th Hodge--Tate twist
$\R(n)$, for $n\in \Z$.
Similarly, we define the dynamical homology $H_1^{dyn}$ as the
graded vector space given by the sum of the terms in the
filtration of $H_1({\mathcal S}_T)$, introduced in Proposition
\ref{homology-dyn}. These vector spaces are again considered with
twisted $\R(n)$-coefficients. The pair

\begin{equation}\label{hom-and-cohom}
 H^1_{dyn} \oplus H_1^{dyn}
\end{equation}
provides a geometric setting, defined in terms of the dynamics of
the shift operator $T$, which contains a copy of the Archimedean
cohomology of \cite{KC}  and of its dual.

\begin{defn}
Let $H^1({\mathcal S}_T,\kappa) =\lim_n F_n$, for a filtration $F_n$
as in Theorem \ref{thm-homologydyn}, with real or complex
coefficients. Let ${\rm Gr}_n =F_n/F_{n-1}$ be the corresponding
graded pieces, with ${\rm Gr}_0 = F_0$.
\begin{enumerate}
\item We define a graded linear subspace ${\mathcal V}$ of the Hilbert
space ${\mathcal L}$, as the span of the elements
$$ \hat\Pi_n \chi_{{\mathcal S}^+(w_{n,k})}, $$
with $\chi_{{\mathcal S}^+(w_{n,k})}$ the characteristic function of
${\mathcal S}^+(w_{n,k})\subset {\mathcal S}^+$ with
$w_{n,k}:=a_0 a_1\ldots a_{n-1}=\underbrace{g_k g_k \ldots
g_k}_{n-times}$. The operator $\hat\Pi_n$ is the  projection
$\hat\Pi_n =\Pi_n- \Pi_{n-1}$, with $\Pi_n$ the orthogonal projection
of ${\mathcal L}$ onto ${\mathcal P}_n$.
\item  We define the dynamical cohomology as
\begin{equation} \label{H-dyn}
H^1_{dyn} := \oplus_{n\leq 0} gr_{2n}^\Gamma H^1_{dyn},
\end{equation}
where we set
\begin{equation} \label{H-dyn-p}
gr_{2n}^\Gamma H^1_{dyn} := {\rm Gr}_{-n} \otimes_{\R} \R(n)
\end{equation}
with $\R(n)= (2\pi \sqrt{-1})^n \R$.

Furthermore, we define the graded subspace of $H^1_{dyn}$
\[ \bar{\mathcal V }:= \oplus_{n\le 0}gr^\Gamma_{2n}\bar{\mathcal V} \]
where $gr^\Gamma_{2n}\bar{\mathcal V}$ is generated by the elements $(2\pi
\sqrt{-1})^n \chi_{-n+1,k}$, for $\chi_{n,k}:=[ \chi_{{\mathcal
S}^+(w_{n,k})} ]\in {\rm Gr}_{n-1}$. \label{def-H-dyn}

\item The dynamical homology $H_1^{dyn}$ is defined as
\begin{equation}\label{H1dyn-hom}
H_1^{dyn}:= \oplus_{n\geq 1} gr_{2n}^\Gamma H_1^{dyn},
\end{equation}
where we set
\begin{equation}\label{H-hom-dyn-gr}
gr_{2n}^\Gamma H_1^{dyn} := {\mathcal K}_{n-1}\otimes_\R \R(n).
\end{equation}
We also define $\mathcal W \subset H_1^{dyn}$ as the graded
sub-space $\mathcal W = \oplus_{n\ge 1}gr^\Gamma_{2n}\mathcal W$,
where $gr_{2n}^\Gamma\mathcal W$ is generated by the $2g$ elements
$(2\pi\sqrt{-1})^n \, \, \underbrace{g_k g_k \ldots
g_k}_{n-times}$.
\end{enumerate}
\end{defn}

The choice of indexing the grading by $gr^\Gamma_{2n}$ instead of
$gr^\Gamma_{n}$ is motivated by comparison to the grading on the
cohomological construction of \cite{KC}.

In \cite{CM}, we showed that the subspaces ${\mathcal V}$ and
$\bar{\mathcal V}$ realize copies of the Archimedean cohomology of
\cite{KC} embedded in the space of cochains ${\mathcal L}$ of the
dynamical cohomology and in the dynamical cohomology itself, while
the pair $\bar{\mathcal V}\oplus {\mathcal W}$ realize a copy of
the cohomology of the cone of the ``local monodromy map'' $N$ of
\cite{KC} inside the pair of dynamical cohomology and homology
$H^1_{dyn}\oplus H_1^{dyn}$. The isomorphism between ${\mathcal
V}$ and the Archimedean cohomology is realized by a natural choice
of a basis of holomorphic differentials for the Archimedean
cohomology, constructed from the data of the Schottky
uniformozation as in \cite{Man}.

An explicit geometric description for the space ${\mathcal V}$ in
terms of geodesics in $\mX_\Gamma$ is obtained by interpreting the
characteristic function $\chi_{{\mathcal S}^+(w_{n,k})}$ as the
``best approximation'' within ${\mathcal P}_n$ to a distribution
supported on the periodic sequence of period $g_k$,
$$ f(g_k g_k g_k g_k g_k \ldots) =1 \ \ \text{ and } \ \ f(a_0 a_1 a_2 a_3
\ldots)=0 \, \, \text{ otherwise.} $$ Such periodic sequence $g_k
g_k g_k g_k g_k \ldots$ describes the closed geodesic in
$\mX_\Gamma$ that is the oriented core of one of the handles in
the handlebody. Thus, the  subspace $gr_{2n}^\Gamma {\mathcal V}
\subset {\mathcal P}_n$ is spanned by the ``best
approximations'' within ${\mathcal P}_n$ to cohomology
classes supported on the core handles of the handlebody. In other
words, this interpretation views the index $n\leq 0$ of the graded
structure ${\mathcal V}=\oplus_n gr_{2n}^\Gamma {\mathcal V}$ as a
measure of ``zooming in'', with increasing precision for larger
$|n|$, on the core handles of the handlebody $\mX_\Gamma$.

\subsection{A spectral triple from dynamics.}

Recall that a spectral triple consists of the following data (\cf
\cite{Connes}).

\begin{defn}\label{specDef}
a spectral triple $({\mathcal A}, {\mathcal H}, D)$ consists of
a ${\rm C}^*$-algebra ${\mathcal A}$ with a representation
$$ \rho : {\mathcal A} \to {\mathcal B}({\mathcal H}) $$
as bounded operators on a Hilbert space ${\mathcal H}$, and an
operator  $D$ (called the Dirac operator) on ${\mathcal H}$, which
satisfies the following properties:
\begin{enumerate}
\item $D$ is self--adjoint.
\item For all $\lambda\notin \R$, the resolvent $(D-\lambda)^{-1}$ is
a compact operator on ${\mathcal H}$.
\item The commutator $[D,a]$ is a bounded
operator on ${\mathcal H}$, for all $a\in {\mathcal A}_0$, a dense
involutive subalgebra of ${\mathcal A}$.
\end{enumerate}
\end{defn}

We consider the Cuntz--Krieger algebra ${\mathcal O}_A$ (\cf \cite{Cu}
\cite{CuKrie}) defined as the universal ${\rm C}^*$--algebra
generated by partial isometries $S_1, \ldots, S_{2g}$, satisfying
the relations
\begin{equation} \label{CK1rel} \sum_j S_j S_j^* =I \end{equation}
\begin{equation} \label{CK2rel} S_i^* S_i =\sum_j A_{ij} \, S_j
S_j^*, \end{equation} where $A=(A_{ij})$ is the $2g\times 2g$
transition matrix of the subshift of finite type $({\mathcal
S},T)$, namely the matrix whose entries are $A_{ij}=1$ whenever
$|i-j|\neq g$, and $A_{ij}=0$ otherwise.

The algebra ${\mathcal O}_A$ can be also described in terms of the
action of the free group $\Gamma$ on its limit set
$\Lambda_\Gamma$ (\cf \cite{Rob}, \cite{Spi}), so that we can
regard ${\mathcal O}_A$ as a noncommutative space replacing the
classical quotient $\Lambda_\Gamma / \Gamma$. In fact,
the action of $\Gamma$ on $\Lambda_\Gamma \subset \P^1(\C)$
determines a unitary representation of ${\mathcal O}_A$ on the Hilbert
space $L^2(\Lambda_\Gamma, \mu)$, given by
\begin{equation} (T_{\gamma^{-1}} f)(x)
:=|\gamma^\prime(x)|^{\delta_H/2}\, f(\gamma x),  \ \ \ \text{ and } \ \ \
(P_\gamma f)(x) := \chi_\gamma (x) f(x), \label{TandP}
\end{equation}
where $\delta_H$ is the Hausdorff dimension of $\Lambda_\Gamma$ and the
element $\gamma\in \Gamma$ is identified with a reduced word in the
generators $\{ g_j \}_{j=1}^g$ and their inverses, and
$\chi_\gamma$ is the characteristic function of the cylinder
$\Lambda_\Gamma(\gamma)$ of all (right) infinite reduced words
that begin with the word $\gamma$.
This determines an identification of ${\mathcal O}_A$ with the
(reduced) crossed product ${\rm C}^*$--algebra, ${\mathcal O}_A
\cong C(\Lambda_\Gamma) \rtimes \Gamma$.

We then consider on the Hilbert space ${\mathcal
L}=L^2(\Lambda_\Gamma, \mu)$, the
unbounded linear self adjoint operator
$D: {\mathcal L}\to {\mathcal L}$ given by the grading operator of
the filtration ${\mathcal P}_n$, namely,
\begin{equation}\label{D-L}
D = \sum_n n\, \hat\Pi_n.
\end{equation}

The restriction of this operator to the subspace ${\mathcal V}$ of
the dynamical cohomology, isomorphic to the Archimedean cohomology
of \cite{KC}, agrees with the ``Frobenius'' operator $\Phi$
considered in \cite{KC}, which computes the local factor as a
regularized determinant as in \cite{Den}.

We extend the operator \eqref{D-L} to an operator ${\mathcal D}$ on
${\mathcal H}={\mathcal L}\oplus {\mathcal L}$ as
\begin{equation}\label{D-oper3L}
{\mathcal D}|_{{\mathcal L}\oplus 0} = \sum_n (n+1) \, (\hat\Pi_n \oplus
0)  \ \ \ \
{\mathcal D}|_{0\oplus {\mathcal L}} = -\sum_n n \,\, (0\oplus\hat\Pi_n).
\end{equation}

The presence of a shift by one in the grading operator,
reflects the shift by one
in the grading that appears in the duality isomorphisms on the
cohomology of the cone of the monodromy map $N$ as in Proposition 4.8
of \cite{KC}. The Dirac operator
\eqref{D-oper3L} takes into account the presence of this shift.

\smallskip

There is another possible natural choice for the sign of the Dirac
operator, instead of the one in \eqref{D-oper3L}. Instead of being
determined by the sign of the operator $\Phi$ on the archimedean
cohomology and its dual, this other choice is determined by the
duality map that exists on the complex of 
\cite{KC}, realized by powers of the monodromy map (\cf \cite{KC},
Proposition 4.8). In this case, the sign would then be given by the
operator 
$$ F = \left(\begin{array}{cc} 0 & 1 \\ 1 & 0 \end{array}\right), $$
that exchanges the two copies of ${\mathcal L}$ in ${\mathcal H}$.

\medskip

We assume that the Schottky group $\Gamma$ has limit set
$\Lambda_\Gamma$ of Hausdorff dimension $\delta_H <1$. Then
the data above define the ``dynamical spectral triple'' at arithmetic
infinity.

\begin{thm}\label{dyn-SP3OA}
The data $({\mathcal O}_A,{\mathcal H},{\mathcal D})$, where the
algebra ${\mathcal O}_A$ acts diagonally on ${\mathcal H}={\mathcal
L}\oplus {\mathcal L}$, and the operator ${\mathcal D}$ is given by
\eqref{D-oper3L} form a spectral triple in the
sense of Connes, as in Definition \ref{specDef}.
\end{thm}

The bound on the commutators with the generators $S_i$ of ${\mathcal
O}_A$ and their adjoints, is obtained in \cite{CM} in terms of the
Poincar\'e series of the Schottky group.

The spectral triple defined this way appears to be related to spectral
triples for AF algebras, recently introduced in \cite{AntChris}. In such
constructions the Dirac operator generalizes the grading operator
\eqref{D-L}, by operators of the form $D=\sum \alpha_n \hat\Pi_n$,
where $\hat\Pi_n=\Pi_n - \Pi_{n-1}$ are the projections associated to
a filtration of the AF algebra, and the coefficients $\alpha_n$ given
by a sequence of positive real numbers, satisfying certain growth
conditions.

The $C^*$--algebra ${\rm C}(\Lambda_\Gamma)$ is a commutative
AF--algebra (approximately finite dimensional), obtained as the direct
limit of the finite dimensional commutative $C^*$--algebras generated
by characteristic functions of a covering of $\Lambda_\Gamma$. This
gives rise to the filtration ${\mathcal P}_n$ in Theorem
\ref{thm-homologydyn}, hence our choice of Dirac operator fits into
the setting of \cite{AntChris} for the AF algebra ${\rm
C}(\Lambda_\Gamma)$. On the other hand, while in the construction of
\cite{AntChris} the eigenvalues $\alpha_n$ can be chosen sufficiently
large, so that the resulting spectral triple for the AF algebra would
be finitely summable, when we consider the Cuntz--Krieger algebra
${\mathcal O}_A \cong {\rm C}(\Lambda_\Gamma)\rtimes \Gamma$, the
boundedness of commutators (condition {\em 3.} of Definition
\ref{specDef}) can only be satisfied for the special choice of
$\alpha_n = c \, n$, with $c$ a constant, \cf Remark 2.2 of
\cite{AntChris}, which does not yield a finitely summable spectral
triple. The reason for this lies in a well known result of Connes
\cite{Connes2} which shows that non amenable discrete groups (as is
the case for the Schottky group $\Gamma$) do not admit finitely
summable spectral triples. Thus, if the dense subalgebra of ${\mathcal
O}_A$ with which ${\mathcal D}$ has bounded commutators contains group
elements, then the Dirac operator ${\mathcal D}$ cannot be finitely
summable. 

In our construction, the choice of
sign for ${\mathcal D}$ is prescribed by the graded structure of
the cohomology theory of \cite{KC} and by the identification of the
Archimedean cohomology of \cite{KC} and its dual with subspaces of
the dynamical cohomology and homology as in \cite{CM}. This way,
the requirement that the Dirac operator agrees with the operator
$\Phi$ of \cite{KC} on these subspaces fixes the choice of the
sign of the Dirac operator, which carries the topological
information on the noncommutative manifold. 
In the construction of spectral triples for AF algebras
of \cite{AntChris} only the metric aspect of the spectral triple
is retained, that is, the operator considered is of the form
$|D|$, while the sign is not discussed.

A possible way to refine the construction of the spectral triple and
deal with the lack of finite summability is 
through the fact that the Cuntz--Krieger
algebra ${\mathcal O}_A$ has a second description as a crossed product
algebra. Namely, up to stabilization (\ie tensoring with compact
operators) we have
\begin{equation} \label{AF-T} {\mathcal O}_A \simeq
{\mathcal F}_A \rtimes_T \Z, \end{equation}
where ${\mathcal F}_A$ is an approximately finite
dimensional (AF) algebra, \cf \cite{Cu},
\cite{CuKrie}. This algebra can be described in terms of a groupoid
${\rm C}^*$-algebra associated to the ``unstable manifold'' in the
Smale space $({\mathcal S},T)$. In fact, consider the algebra
${\mathcal O}_A^{alg}$ generated algebraically by the $S_i$ and
$S_i^*$ subject to the Cuntz--Krieger relations \eqref{CK1rel}
\eqref{CK2rel}. Elements in ${\mathcal O}_A^{alg}$ are linear
combinations of monomials $S_\mu S_\nu^*$, for multi-indices $\mu$,
$\nu$, \cf \cite{CuKrie}. The AF algebra ${\mathcal F}_A$ is generated
by elements $S_\mu S_\nu^*$ with $|\mu|=|\nu|$, and is filtered by
finite dimensional algebras ${\mathcal F}_{A,n}$ generated by elements
of the form $S_\mu P_i S_\nu^*$ with $|\mu|=|\nu|=n$ and
$P_i=S_iS_i^*$ the range projections, and embeddings determined by the
matrix $A$. The commutative algebra ${\rm C}(\Lambda_\Gamma)$ sits as
a subalgebra of ${\mathcal F}_A$ generated by all range projections
$S_\mu S_\mu^*$. The embedding is compatible with the filtration and
with the action of the shift $T$, which is implemented on ${\mathcal
F}_A$ by the transformation $a\mapsto \sum_i S_i\, a\, S_i^*$. (\cf
\cite{CuKrie}.)

The fact that the algebra can be written in the form \eqref{AF-T}
implies that, by Connes' result on hyperfiniteness \cite{Connes2},
it may carry a finitely summable spectral triple. It
is an interesting problem whether the construction of a finitely
summable triple can be carried out in a way that is of arithmetic
significance.

\subsection{Local factor}

For an arithmetic variety ${\mathcal X}$ over $\Sp \Z$, the
``Archimedean factor'' (local factor at arithmetic infinity)
$L_\kappa(H^m,s)$ is a product of Gamma functions, with exponents and
arguments that depend on the Hodge structure
$H^m = H^m (X,\C)= \oplus_{p+q=m} H^{p,q}$. More precisely,
(\cf \cite{Serre})
\begin{equation}\label{factor}
 L_\kappa (H^m,s)= \begin{cases}
\prod_{p,q}\Gamma_\C(s-\text{min}(p,q))^{h^{p,q}}
&\kappa = \C  \\\\
\prod_{p<q}\Gamma_\C(s-p)^{h^{p,q}}\prod_p
\Gamma_\R(s-p)^{h^{p+}}\Gamma_\R(s-p+1)^{h^{p-}}
& \kappa = \R,
\end{cases}
\end{equation}
where the $h^{p,q}$, with $p+q=m$, are the Hodge numbers,
$h^{p,\pm}$ is the dimension of the $\pm(-1)^p$-eigenspace
of de Rham conjugation on $H^{p,p}$, and
$$ \Gamma_\C(s) := (2\pi)^{-s}\Gamma(s) \ \ \ \ \
\Gamma_\R(s) :=2^{-1/2}\pi^{-s/2}\Gamma(s/2). $$

Deninger produced a unified description of the factors at arithmetic
infinity and at the finite primes, in the form of a Ray--Singer
determinant (\cf \cite{Den}). The factor \eqref{factor} satisfies
\begin{equation}\label{DenDet}
L_\kappa(H^m,s) =
\det_{\infty}\left(\frac{1}{2\pi}(s-\Phi)|_{{\mathcal V}^m}
\right)^{-1},
\end{equation}
where ${\mathcal V}^m$ is an infinite dimensional vector space. The zeta
regularized determinant of an unbounded self adjoint operator $T$ is
defined as $\det_\infty(s-T)=\exp (-\frac{d}{dz}
\zeta_T(s,z)|_{z=0})$.

In \cite{KC}, the spaces ${\mathcal V}^\cdot$ are identified with
inertia invariants of a double complex of real Tate-twisted
differential forms on $X_{/\C}$ with suitable cutoffs. Namely,
such complex is endowed with the action of an endomorphism  $N$,
which represents a ``logarithm of the local monodromy at
arithmetic infinity'', and the spaces ${\mathcal V}^\cdot$ are
identified with the kernel of the map $N$ on the hypercohomology.
In particular, in the case of an arithmetic surface, we showed in
\cite{CM} that the Archimedean cohomology group ${\mathcal V}^1$
is identified with a subspace ${\mathcal V}$ of ${\mathcal P}$ and
of the dynamical cohomology $H^1_{dyn}$.

The dynamical spectral triple of Theorem
\ref{dyn-SP3OA} is not finitely summable.
However, it is possible to recover from these data the local
factor at arithmetic infinity \eqref{DenDet} for $m=1$.

\begin{prop}\label{L-factor2}
Consider the zeta functions
\begin{equation}\label{zeta-proj1}
\zeta_{\pi({\mathcal V}), D} (s,z):=
\sum_{\lambda \in \Sp(D)} \Tr \left( \pi({\mathcal V})
\Pi(\lambda, D) \right) (s-\lambda)^{-z},
\end{equation}
for $\pi({\mathcal V})$ the orthogonal projection on
the norm closure of ${\mathcal V}$ in ${\mathcal L}$, and
\begin{equation}\label{zeta-proj2}
\zeta_{\pi({\mathcal V}, \bar F_\infty=id), D} (s,z):=
\sum_{\lambda \in \Sp(D)} \Tr \left( \pi({\mathcal V}, \bar
F_\infty=id) \Pi(\lambda, D) \right) (s-\lambda)^{-z},
\end{equation}
for $\pi({\mathcal V}, \bar F_\infty=id)$ the orthogonal
projection on the norm closure of ${\mathcal V}^{\bar
F_\infty=id}$. The corresponding regularized
determinants satisfy
\begin{equation}\label{L-det21}
\exp\left( - \frac{d}{dz} \zeta_{\pi({\mathcal V}), D/2\pi}
(s/2\pi,z)|_{z=0} \right)^{-1} =
L_\C(H^1(X),s),
\end{equation}
\begin{equation}\label{L-det22}
\exp\left( - \frac{d}{dz} \zeta_{\pi({\mathcal V}, \bar
F_\infty=id), D/2\pi} (s/2\pi,z)|_{z=0} \right)^{-1} =
L_\R(H^1(X),s).
\end{equation}
Moreover, the operator $\pi({\mathcal V})$ acts
on the range of the spectral projections
$\Pi(\lambda, D)$ as certain elements of the
algebra ${\mathcal O}_A$.
\end{prop}

Here $\bar F_\infty$ is the involution induced by the real structure
on $X_{/\R}$, which corresponds to the change of orientation on the
geodesics in $\mX_\Gamma$ and on ${\mathcal S}_T$.

\section{A Dynamical theory for Mumford curves.}

Throughtout this chapter $K$ will denote a finite extension of
$\Q_p$ and $\Delta_K$ the Bruhat-Tits tree associated to $G =
\PGL(2,K)$. In the following we recall few results about the
action of a Schottky group on a Bruhat-Tits tree and on ${\rm
C}^*$-algebras of graphs. Detailed explanations are contained in
\cite{Ma}, \cite{Mum} and \cite{BPRS}, \cite{KuPa},
\cite{KuPaRae}, \cite{Spi}.

Recall that the Bruhat--Tits tree is constructed as follows. One
considers the set of free $\O$-modules of rank $2$:
$M\subset V$. Two such modules are {\it equivalent} $M_1 \sim M_2$ if
there exists an element $\lambda\in K^\ast$, such that $M_1 = \lambda
M_2$.
The group $\rm{GL}(V)$ of linear automorphisms of $V$ operates on the
set of such modules {\it on the left}: $gM = \{gm~|~m\in M\}$, $g\in
\rm{GL}(V)$. Notice that the relation $M_1 \sim M_2$ is equivalent to
the condition that $M_1$ and $M_2$ belong to the same orbit of the center
$K^\ast \subset \rm{GL}(V)$. Hence, the group $G = \rm{GL}(V)/K^\ast$
operates (on the left) on the set of classes of equivalent modules.

We denote by $\Delta^0_K$ the set of such classes and by $\{M\}$
the class of the module $M$. Because $\O$ is a principal ideals
domain and every module $M$ has two generators, it follows that
\[
\{M_1\}, \{M_2\} \in\Delta^0_K, M_1\supset
M_2\quad\Rightarrow\quad M_1/M_2 \simeq \O/\m^l \oplus
\O/\m^k,\quad l,k\in\N.
\]

The multiplication of $M_1$ and $M_2$ by elements of $K$ preserves the
inclusion $M_1\supset M_2$, hence the natural number
\begin{equation}\label{dist}
d(\{M_1\},\{M_2\}) = \vert l-k\vert
\end{equation}
is well defined.

The graph $\Delta_K$ of the group $\PGL(2,K)$ is the
infinite graph with set of vertices $\Delta^0_K$, in which two
vertices $\{M_1\},\{M_2\}$ are adjacent and hence connected by an
edge if and only if $d(\{M_1\},\{M_2\}) = 1$. (\cf \cite{Ma} and
\cite{Mum}.)

For a Schottky group $\Gamma \subset \PGL(2,K)$ there is a
smallest subtree $\Delta'_\Gamma \subset \Delta_K$ containing the
axes of all elements of $\Gamma$. The set of ends of
$\Delta'_\Gamma$ in $\P^1(K)$ is $\Lambda_\Gamma$, the limit set
of $\Gamma$. The group $\Gamma$ carries $\Delta'_\Gamma$ into
itself so that the quotient $\Delta'_\Gamma /\Gamma$ is a finite
graph that coincides with the dual graph of the closed fibre of
the minimal smooth model of the algebraic curve $C/K$
holomorphically isomorphic to $X_\Gamma := \Omega_\Gamma /\Gamma$
(\cf~\cite{Mum} p.~163). There is a smallest tree $\Delta_\Gamma$
on which $\Gamma$ acts and such that $\Delta_\Gamma/ \Gamma$ is
the (finite) graph of the specialization of $C$. The curve $C$  is
a $k$-split degenerate, stable curve. When the genus of the fibers
is at least 2 - \ie when the Schottky group has at least $g\geq 2$
generators - the curve $X_\Gamma$ is called a Schottky--Mumford
curve.

\begin{figure}
\begin{center}
\epsfig{file=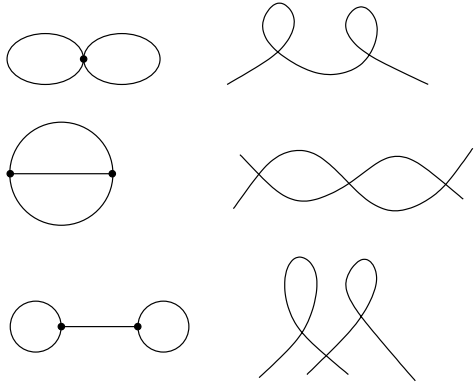} \caption{The graphs
$\Delta_\Gamma/\Gamma$ for genus $g=2$, and the corresponding
fibers. \label{graphs}}
\end{center}
\end{figure}

The possible graphs $\Delta_\Gamma/\Gamma$ and the corresponding
fiber for the case of genus 2 are illustrated in Figure
\ref{graphs}.\smallskip

To a locally finite directed graph one associates a ${\rm
C}^*$-algebra in the following way.

A Cuntz--Krieger family consists of a collection $\{ P_v \}_{v\in
E^0}$ of mutually orthogonal projections and a $\{ S_w \}_{w \in
E^1_+ }$ of partial isometries, satisfying the conditions: $S_w^*
S_w =P_{r(w)}$ and, for all $v\in s(E^1_+)$, $P_v = \sum_{w:
s(w)=v} S_w S_w^*$.

The Cuntz--Krieger elements $\{ P_v, S_w \}$ satisfy the relation
$S_w^* S_w = \sum A_+(w,\tilde w) S_{\tilde w} S_{\tilde w}^*$,
with $A_+$ the edge matrix of the graph.

One defines a universal $C^*$-algebra $C^*(E)$ generated by a
Cuntz--Krieger family. If $E$ is a finite graph with no sinks, we
have $C^*(E)\simeq {\mathcal O}_{A_+}$, where ${\mathcal O}_{A_+}$
is the Cuntz-Krieger algebra of the edge matrix $A_+$. If the
directed graph is a tree $\Delta$, then $C^*(\Delta)$ is an  AF
algebra strongly Morita equivalent to the commutative
$C^*$-algebra $C_0(\partial \Delta)$. A monomorphism of directed
trees induces an injective $*$-homomorphism of the corresponding
$C^*$-algebras.

If $G\subset Aut(E)$ is a group acting freely on the directed
graph $E$, with quotient graph $E/G$, then the crossed product
$C^*$-algebra $C^*(E)\rtimes G$ is strongly Morita equivalent to
$C^*(E/G)$. In particular, if $\Delta$ is the universal covering
tree of a directed graph $E$ and $\Gamma$ is the fundamental
group, then the algebra $C^*(E)$ is strongly Morita equivalent to
$C_0(\partial \Delta)\rtimes \Gamma$.\smallskip

In the following we consider the Bruhat--Tits tree $\Delta_K$ for
a fixed finite extension $K$ of $\Q_p$, and the corresponding
$C^*$-algebra $C^*(\Delta_K)$, which is strongly Morita equivalent
to the abelian $C^*$-algebra of complex valued functions
$C(\P^1(K))$.

Notice that the commutative $C^*$-algebra $C^*(X_K)$ of complex
valued continuous functions on the Mumford curve
$X_K=\Omega_\Gamma /\Gamma$ is strongly Morita equivalent to the
crossed product $C_0(\Omega_\Gamma) \rtimes \Gamma$, with
$C_0(\Omega_\Gamma)$ a $\Gamma$-invariant ideal of $C(\P^1(K))$
with quotient algebra $C(\Lambda_\Gamma)$. The algebra
$C^*(\Delta_\Gamma /\Gamma)$, in turn, is strongly Morita
equivalent to $C^*(\Delta_\Gamma)\rtimes \Gamma$ and to
$C(\Lambda_\Gamma)\rtimes \Gamma$, where $\partial \Delta_\Gamma =
\Lambda_\Gamma$. Similarly, one sees that the algebra
$C^*(\Delta_K/\Gamma)$ is strongly Morita equivalent to the
crossed product algebra $C^*(\Delta_K)\rtimes \Gamma$, which in
turn is strongly Morita equivalent to $C^*(\P^1(K))\rtimes
\Gamma$.

Thus, up to Morita equivalence, the graph algebra
$C^*(\Delta_K/\Gamma)$ can be regarded as a way of extending
the commutative $C^*$-algebra $C^*(X_K)$ (functions on the Mumford
curve) by the Cuntz-Krieger algebra $C^*(\Delta_\Gamma
/\Gamma)$ associated to the edge matrix of the finite graph
$\Delta_\Gamma /\Gamma$.\smallskip

In this paragraph we introduce a dynamical system associated to
the space $\cW(\Delta/\Gamma)$ of walks on the directed tree
$\Delta$ on which $\Gamma$ acts. In particular, we are interested
in the cases when $\Delta = \Delta_K,~\Delta_\Gamma$.\smallskip

For $\Delta=\Delta_\Gamma$, we obtain a subshift of finite type
associated to the action of the Schottky group $\Gamma$ on the
limit set $\Lambda_\Gamma$, of the type that was considered in
\cite{CM}.

Let $\bar V \subset \Delta_\Gamma$ be a finite subtree whose set
of edges consists of one representative for each $\Gamma$-class.
This is a fundamental domain for $\Gamma$ in the weak sense
(following the notation of \cite{Ma}), since some vertices may be
identified under the action of $\Gamma$. Correspondingly,
$V\subset \P^1(K)$ is the set of ends of all infinite paths
starting at points in $\bar V$.

Consider the set $\cW(\Delta_\Gamma/\Gamma)$ of doubly infinite
walks on the finite graph $\Delta_\Gamma/\Gamma$. These are doubly
infinite admissible sequences in the finite alphabet given by the
edges of $\bar V$ with both possible orientations. On
$\cW(\Delta_\Gamma/\Gamma)$ we consider the topology generated by
the sets $\cW^s(\omega,\ell)=\{ \tilde\omega \in
\cW(\Delta_\Gamma/\Gamma): \tilde\omega_k=\omega_k, k\geq \ell\}$
and $\cW^u(\omega,\ell)=\{ \tilde\omega \in
\cW(\Delta_\Gamma/\Gamma): \tilde\omega_k=\omega_k, k\leq \ell\}$,
for $\omega\in \cW(\Delta_\Gamma/\Gamma)$ and $\ell \in \Z$. With
this topology, the space $\cW(\Delta_\Gamma/\Gamma)$ is a totally
disconnected compact Hausdorff space.

The invertible shift map $T$, given by $(T\omega)_k =
\omega_{k+1}$, is a homeomorphism of $\cW(\Delta_\Gamma/\Gamma)$.
We can describe again the dynamical system
$(\cW(\Delta_\Gamma/\Gamma),T)$ in terms of subshifts of finite
type.

\begin{lem}\label{subshift}
The space $\cW(\Delta_\Gamma/\Gamma)$ with the action of the
invertible shift $T$ is a subshift of finite type, where
$\cW(\Delta_\Gamma/\Gamma)=\cS_A$ with $A$ the directed edge
matrix of the finite graph $\Delta_\Gamma /\Gamma$.
\end{lem}

We consider the mapping torus of $T$:
\begin{equation}\label{map-torus-G}
\cW(\Delta_\Gamma/\Gamma)_{T}:= \cW(\Delta_\Gamma/\Gamma) \times
[0,1] / (Tx,0)\sim (x, 1).
\end{equation}\medskip

\subsection{Genus two example}

\begin{figure}
\begin{center}
\epsfig{file=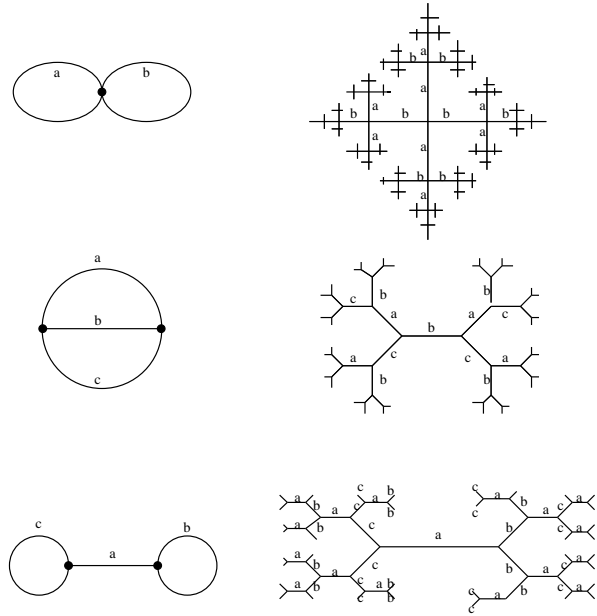} \caption{The graphs $\Delta_\Gamma/\Gamma$
for genus $g=2$, and the corresponding trees $\Delta_\Gamma$.
\label{trees}}
\end{center}
\end{figure}

In the example of Mumford--Schottky curves of genus $g=2$, the
tree $\Delta_\Gamma$ is illustrated in Figure \ref{trees}.

In the first case in the Figure, the tree $\Delta_\Gamma$ is just
a copy of the Cayley graph of the free group $\Gamma$ on two
generators, hence we can identify doubly infinite walks in
$\Delta_\Gamma$ with doubly infinite reduced words in the
generators of $\Gamma$ and their inverses. The directed edge
matrix is given by
$$ A= \left( \begin{array}{cccc} 1 &1& 0& 1\\
1&1&1& 0  \\
0&1&1& 1 \\
1& 0& 1&1  \end{array} \right). $$

In the second case in Figure \ref{trees}, we label by $a=e_1$,
$b=e_2$ and $c=e_3$ the oriented edges in the graph
$\Delta_\Gamma/\Gamma$, so that we have a corresponding set of
labels $E=\{ a,b,c,\bar a, \bar b, \bar c \}$ for the edges in the
covering $\Delta_\Gamma$. A choice of generators for the group
$\Gamma \simeq \Z * \Z$ acting on $\Delta_\Gamma$ is obtained by
identifying the generators $g_1$ and $g_2$ of $\Gamma$ with the
chains of edges $a \bar b$ and $a \bar c$. Doubly infinite walks
in the tree $\Delta_\Gamma$ are admissible doubly infinite
sequences of such labels, where admissibility is determined by the
directed edge matrix
$$ A= \left( \begin{array}{cccccc} 0 &1& 0& 0& 0& 1 \\
1&0&1& 0 & 0 & 0 \\
0&1&0& 1& 0& 0 \\
0& 0& 1& 0&1&0 \\
0&0&0&1&0&1 \\
1&0&0&0&1&0 \end{array} \right). $$

The third case in Figure \ref{trees} is analogous. A choice of
generators for the group $\Gamma\simeq \Z * \Z$ acting on
$\Delta_\Gamma$ is given by $ab\bar a$ and $c$. Doubly infinite
walks in the tree $\Delta_\Gamma$ are admissible doubly infinite
sequences in the alphabet $E=\{ a,b,c,\bar a, \bar b, \bar c \}$,
with admissibility determined by the directed edge matrix
$$ A= \left( \begin{array}{cccccc}
0 & 0 & 1 & 0 & 0 & 1 \\
1 & 1 & 0 & 0 & 0 & 0 \\
0 & 0 & 1 & 1 & 0 & 0 \\
0 & 1 & 0 & 0 & 1 & 0 \\
1 & 0 & 0 & 0 & 1 & 0 \\
0 & 0 & 0 & 1 & 0 & 1 \end{array} \right). $$

The construction is analogous for genus $g>2$, for the various
possible finite graphs $\Delta_\Gamma/\Gamma$. The directed edge
matrix can then be written in block form as
$$ A=\left( \begin{array}{cc} \alpha_{11} & \alpha_{12} \\ \alpha_{21}
& \alpha_{22} \end{array} \right), $$ where each block
$\alpha_{ij}$ is a $\# (\Delta_\Gamma /\Gamma)^1_+ \times \#
(\Delta_\Gamma /\Gamma)^1_+$--matrix with
$\alpha_{12}=\alpha_{12}^t$, $\alpha_{21}=\alpha_{21}^t$, and
$\alpha_{11}= \alpha_{22}^t$.

\subsection{Cohomology of $\cW(\Delta/\Gamma)_{T}$.}\label{A}

Let $\Delta=\Delta_\Gamma$. We identify the first cohomology group
$H^1(\cW(\Delta_\Gamma/\Gamma)_{T},\Z)$ with the group of homotopy
classes of continuous maps of $\cW(\Delta_\Gamma/\Gamma)_{T}$ to
the circle. Let ${\rm C}(\cW(\Delta_\Gamma/\Gamma),\Z)$ be the
$\Z$-module of integer valued continuous functions on
$\cW(\Delta_\Gamma/\Gamma)$, and let
$$
{\rm C}(\cW(\Delta_\Gamma/\Gamma),\Z)_T:= Coker(\delta),
$$
for $\delta(f)=f-f\circ T$. The analog of Theorem
\ref{thm-homologydyn} holds:

\begin{prop}\label{filtration}
The map $f\mapsto [\exp(2\pi i t f(x))]$, which associates to an
element $f\in {\rm C}(\cW(\Delta_\Gamma/\Gamma),\Z)$ a homotopy
class of maps from $\cW(\Delta_\Gamma/\Gamma)_{T}$ to the circle,
gives an isomorphism ${\rm C}(\cW(\Delta_\Gamma/\Gamma),\Z)_T
\simeq H^1(\cW(\Delta_\Gamma/\Gamma)_{T},\Z)$. Moreover, there is
a filtration of ${\rm C}(\cW(\Delta_\Gamma/\Gamma),\Z)_T$ by free
$\Z$-modules $F_0\subset F_1\subset \cdots F_n \cdots$, of rank
$\theta_n - \theta_{n-1} +1 $, where $\theta_n$ is the number of
admissible words of length $n+1$ in the alphabet, so that we have
$$ H^1(\cW(\Delta_\Gamma/\Gamma)_{T},\Z)=\varinjlim_n
F_n. $$ The quotients $F_{n+1}/F_n$ are also torsion free.
\end{prop}

The space $\cW(\Delta_\Gamma/\Gamma)_T$ corresponds to a space of
``bounded geodesics'' on the graph $\Delta_K/\Gamma$, where
geodesics, in this setting, are just doubly infinite walks in
$\Delta_K/\Gamma$. In particular, a closed geodesic is the image
under the quotient map $\pi_\Gamma: \Delta_K \to \Delta_K/\Gamma$
of a doubly infinite walk in the Bruhat-Tits tree $\Delta_K$ with
ends given by the pair $z^+(\gamma), z^-(\gamma)$ of fixed points
of some element $\gamma \in \Gamma$.  Similarly, a bounded
geodesic is an element $\omega \in \cW(\Delta_K/\Gamma)$ which is
the image, under the quotient map, of a doubly infinite walk in
$\Delta_K$ with both ends on $\Lambda_\Gamma \subset \P^1(K)$.
This implies that a bounded geodesic is a walk of the form
$\omega=\pi_\Gamma(\tilde \omega)$, for some $\tilde\omega \in
\cW(\Delta_\Gamma/\Gamma)$. By construction, any such walk is an
axis of $\Delta_\Gamma$.

Orbits of $\cW(\Delta_\Gamma/\Gamma)$ under the action of the
invertible shift $T$ correspond bijectively to orbits of the
complement of the diagonal in $\Lambda_\Gamma \times
\Lambda_\Gamma$ under the action of $\Gamma$. Thus, we see that
$\cW(\Delta_\Gamma/\Gamma)_T$ gives a geometric realization of the
space of ``bounded geodesics'' on the graph $\Delta_K/\Gamma$,
much as, in the case of the geometry at arithmetic infinity, we
used the mapping torus of the shift $T$ as a model of the tangle
of bounded geodesics in a hyperbolic handlebody.

As in the case at infinity, we can consider the Pimsner--Voiculescu
exact sequence computing the $K$-theory groups of the crossed product
${\rm C}^*$-algebra ${\rm C}(\cW(\Delta_\Gamma/\Gamma))\rtimes_T \Z$,
\begin{equation}\label{PVp}
0 \to H^0(\cW(\Delta_\Gamma/\Gamma)_T,\Z) \to {\rm
C}(\cW(\Delta_\Gamma/\Gamma),\Z)
\stackrel{\delta=1-T}{\longrightarrow} {\rm
C}(\cW(\Delta_\Gamma/\Gamma),\Z) \to
H^1(\cW(\Delta_\Gamma/\Gamma)_T,\Z) \to 0.
\end{equation}
In the corresponding sequence
\begin{equation}\label{PVHp}
0 \to H^0(\cW(\Delta_\Gamma/\Gamma)_T,\kappa)   \to {\mathcal P}
\stackrel{\delta}{\longrightarrow} {\mathcal P} \to
H^1(\cW(\Delta_\Gamma/\Gamma)_T,\kappa) \to 0,
\end{equation}
for the cohomology for
$H^*(\cW(\Delta_\Gamma/\Gamma)_T,\kappa)$, with
$\kappa=\R$ or $\C$, we can take the vector space ${\mathcal P}$
obtained, as in the case at infinity, by tensoring with $\kappa$ the
$\Z$-module $\cP\subset {\rm
C}(\cW(\Delta_\Gamma/\Gamma),\Z)$ of functions of future
coordinates where $\cP \simeq {\rm
C}(\cW^+(\Delta_\Gamma/\Gamma),\Z)$. This has a filtration
$\cP=\cup_n \cP_n$, where $\cP_n$ is identified with the submodule
of ${\rm C}(\cW^+(\Delta_\Gamma/\Gamma),\Z)$ generated by
characteristic functions of
$\cW^+(\Delta_\Gamma/\Gamma,\rho)\subset
\cW^+(\Delta_\Gamma/\Gamma)$, where $\rho\in
\cW^*(\Delta_\Gamma/\Gamma)$ is a finite walk $\rho=w_0\cdots w_n$
of length $n+1$, and $\cW^+(\Delta_\Gamma/\Gamma,\rho)$ is the set
of infinite paths $\omega\in \cW^+(\Delta_\Gamma/\Gamma)$, with
$\omega_k =w_k$ for $0\leq k\leq n+1$. This filtration defines the
terms $F_n = \cP_n /\delta\cP_{n-1}$ in the filtration of the
dynamical cohomology of the Mumford curve, as in Proposition
\ref{filtration}.
Again, we will use the same notation in the following for the free
$\Z$-module $\cP_n$ of functions of at most $n+1$ future coordinates
and the vector space obtained by tensoring $\cP_n$ by $\kappa$.

We obtain a Hilbert space completion of the space ${\mathcal P}$ of
cochains in \eqref{PVHp} by considering ${\mathcal
L}=L^2(\Lambda_\Gamma, \mu)$ defined with respect to the measure on
$\Lambda_\Gamma=\partial \Delta_\Gamma$ given by assigning its value
on the clopen set $V(v)$,
given by the ends of all paths in $\Delta_\Gamma$ starting at a vertex
$v$, to be
$$ \mu(V(v))=q^{-d(v)-1}, $$
with $q={\rm card}(\O/\m)$.

\bigskip

In \cite{CM} \S 4, we showed how the mapping torus $\cS_T$ of the
subshift of finite type $(\cS,T)$, associated to the limit set of
the Schottky group, maps surjectively to the tangle of bounded
geodesics inside the hyperbolic handlebody, through a map that
resolves all the points of intersection of different geodesics. In
the case of the Mumford curve, where we replace the real
hyperbolic 3-space by the Bruhat--Tits building $\Delta_K$, the
analog of the surjective map from $\cS_T$ to the tangle of bounded
geodesics is a map from $\cW(\Delta_\Gamma/\Gamma)_T$ to the dual
graph $\Delta_\Gamma/\Gamma$. Here is a description of this
map.\smallskip

As before, we write elements of $\cW(\Delta_\Gamma/\Gamma)$ as
admissible doubly infinite sequences
$$\omega= \ldots w_{i_{-m}}\ldots w_{i_{-1}} w_{i_0} w_{i_1} \ldots
w_{i_n}  \ldots,$$ with the $w_{i_k} =\{ e_{i_k}, \epsilon_{i_k}
\}$ oriented edges on the graph $\Delta_\Gamma/\Gamma$. We
consider each oriented edge $w$ of normalized length one, so that
it can be parameterized as $w(t)=\{ e(t), \epsilon \}$, for $0\leq
t \leq 1$, with $\bar w(t)=\{ e(1-t), -\epsilon \}$. Since
$\omega\in {\mathcal S}_A$ is an admissible sequence of oriented
edges we have $w_{i_k}(1)=w_{i_{k+1}}(0) \in \Delta_\Gamma^{(0)}$.

We consider a map of the covering space $\cW(\Delta_\Gamma/\Gamma)
\times \R$ of $\cW(\Delta_\Gamma/\Gamma)_{T}$ to $|\Delta_\Gamma|$
of the form
\begin{equation}\label{tildeE}
 \tilde E(\omega,\tau) = w_{i_{[\tau]}}(\tau-[\tau]).
\end{equation}
Here $|\Delta_\Gamma|$ denotes the geometric realization of the
graph. By construction, the map $\tilde E$ satisfies $\tilde
E(T\omega,\tau)=\tilde E(\omega,\tau +1)$, hence it descends to a
map $E$ of the quotient
\begin{equation}\label{E}
E : \cW(\Delta_\Gamma)_{T} \to |\Delta_\Gamma|.
\end{equation}
We then obtain a map to  $|\Delta_\Gamma /\Gamma|$, by composing
with the quotient map of the $\Gamma$ action, $\pi_\Gamma:
\Delta_\Gamma \to \Delta_\Gamma/\Gamma$, that is,
\begin{equation}\label{barE}
\bar E:= \pi_\Gamma \circ E : \cW(\Delta_\Gamma)_{T} \to
|\Delta_\Gamma/\Gamma|.
\end{equation}

Thus, we obtain the following.

\begin{prop}
The map $\bar E$ of \eqref{barE} is a continuous surjection from
the mapping torus $\cW(\Delta_\Gamma)_{T}$ to the geometric
realization $|\Delta_\Gamma/\Gamma |$ of the finite graph
$\Delta_\Gamma/\Gamma$.
\end{prop}

The fibers of the map \eqref{barE} are explicitly described as
\begin{equation}\label{Efiber}
 \bar E^{-1}(w_i(t)) = \cW(\Delta_\Gamma) (w_i) \times \{ t \} \cup
\cW(\Delta_\Gamma) (\bar w_i) \times \{ 1-t \},
\end{equation}
where $w_i(t)$, for $t\in [0,1]$ is a parameterized oriented edge
in the graph $\Delta_\Gamma/\Gamma$, and $\cW(\Delta_\Gamma)
(w_i)\subset \cW(\Delta_\Gamma)$ consists of
$$ \cW(\Delta_\Gamma) (w_i)=\{ \omega\in \cW(\Delta_\Gamma): \, \,
w_{i_0}=w_i \}, $$ for $\omega=\ldots w_{i_{-m}}\ldots w_{i_{-1}}
w_{i_0} w_{i_1} \ldots w_{i_n} \ldots$.\medskip

\smallskip

We can also consider the construction described above, where, instead
of using the tree $\Delta_\Gamma$, we use larger $\Gamma$-invariant
trees inside the Bruhat-Tits tree $\Delta$.

The set of doubly infinite walks
$\cW(\Delta_K/\Gamma)$ can be identified with the set of
admissible doubly infinite sequences in the oriented edges of a
fundamental domain for the action of $\Gamma$. Therefore, we
obtain the following identification
\begin{equation}\label{walksK}
 \cW(\Delta_K/\Gamma) \simeq V \times
\cW^*(\Delta_\Gamma/\Gamma)\times V \cup V \times
\cW^+(\Delta_\Gamma/\Gamma) \cup \cW^+(\Delta_\Gamma/\Gamma)
\times V \cup \cW(\Delta_\Gamma/\Gamma),
\end{equation}
where we distinguish between walks that wander off from the finite
graph $\Delta_\Gamma/\Gamma$ along one of the paths leading to an
end in $V$, in one or in both directions, and those that stay
confined within the finite graph $\Delta_\Gamma/\Gamma$. Here we
include in $\cW^*(\Delta_\Gamma/\Gamma)=\cup_n
\cW^n(\Delta_\Gamma/\Gamma)$ also the case
$\cW^0(\Delta_\Gamma/\Gamma)=\emptyset$, where the walk in
$\Delta_K/\Gamma$ does not intersect the finite graph
$\Delta_\Gamma/\Gamma$ at all.

In the topology induced by the $p$-adic norm, $\P^1(K)$ is a
totally disconnected compact Hausdorff space, and so is
$\cW(\Delta_K/\Gamma)$ by \eqref{walksK}. Again, we consider the
invertible shift map $T$, which is a homeomorphism of
$\cW(\Delta_K/\Gamma)$, and we form the mapping torus
\begin{equation}\label{map-torus}
\cW(\Delta_K/\Gamma)_{T}:= \cW(\Delta_K/\Gamma) \times [0,1] /
(Tx,0)\sim (x, 1).
\end{equation}

We obtain the analog of Proposition \ref{filtration}, though in
the case of $\Delta_K$ the $\Z$-modules $F_n$ will not be finitely
generated. On the other hand, we can restrict to ``neighborhoods'' of
the tree $\Delta_\Gamma$ inside $\Delta_K$, which correspond to the
reduction maps modulo powers of the maximal ideal.

In fact, in the theory of Mumford curves, it is important to consider
also the reduction modulo powers $\m^n$ of the maximal ideal $\m
\subset \O_K$, which provides infinitesimal neighborhoods of
order $n$ of the closed fiber.

For each $n\geq 0$, we consider a subgraph $\Delta_{K,n}$ of the
Bruhat-Tits tree $\Delta_K$ defined by setting
$$ \Delta_{K,n}^0 := \{ v\in \Delta_K^0 : \,  d(v,
\Delta_\Gamma ')\leq n \}, $$ with respect to the distance
\eqref{dist}, with $d(v,\Delta_\Gamma ') := \inf \{ d(v,\tilde v):
\, \tilde v \in (\Delta_\Gamma ')^0 \}$, and
$$ \Delta_{K,n}^1 :=\{ w\in \Delta_K^1 : \, s(w), r(w)\in
\Delta_{K,n}^0 \}. $$ Thus, we have $\Delta_{K,0} = \Delta_\Gamma
'$. We have $\Delta_K=\cup_n \Delta_{K,n}$.

For all $n\in \N$, the graph $\Delta_{K,n}$ is invariant under the
action of the Schottky group $\Gamma$ on $\Delta$, and the finite
graph $\Delta_{K,n}/\Gamma$ gives the dual graph of the reduction
$X_K \otimes \O/ \m^{n+1}$. They form a directed family with
inclusions $j_{n,m}: \Delta_{K,n} \hookrightarrow \Delta_{K,m}$,
for all $m\geq n$, with all the inclusions compatible with the
action of $\Gamma$.

In \cite{CM1} we introduce the dynamical cohomology and discuss
spectral geometry for the ${\rm
C}^*$-algebras  ${\rm C}^*(\Delta_{K,n}/\Gamma)\simeq {\rm
C}^*(\Delta_{K,n})\rtimes \Gamma$.

\subsection{Spectral triples and Mumford curves.}

We consider the Hilbert space ${\mathcal H}={\mathcal L}\oplus
{\mathcal L}$ and the operator ${\mathcal D}$
defined as
\begin{equation}\label{D-Lp}
{\mathcal D}|_{{\mathcal L}\oplus 0} = -\frac{2\pi}{R \, \log q} \, \sum \,
(n+1) \, \hat\Pi_n \ \ \ \ {\mathcal D}|_{0\oplus {\mathcal L}} =
\frac{2\pi}{R \, \log q} \, \sum \, n \, \hat\Pi_n,
\end{equation}
where $ \hat\Pi_n = \Pi_n - \Pi_{n-1}$ are the orthogonal projections
associated to the filtration $\cP_n$, the integer $R$ is the length
of all the words representing the generators of $\Gamma$ (this can be
taken to be the same for all generators, possibly after
blowing up a finite number of points on the special fiber, as
explained in \cite{CM1}), and $q={\rm card}(\O/\m)$.

\smallskip

The same argument used in \cite{CM} for the case at arithmetic
infinity adapts to the case of Mumford curves to prove the following
result. (Note: the statement below corrects an unfortunate mistake
that occurred in \S 5.4 of \cite{CM1}.)

\begin{thm}\label{triple1}
Consider the tree $\Delta_\Gamma$ of the $p$-adic Schottky group
acting on $\Delta_K$.
\begin{enumerate}
\item There is a representation of the algebra ${\rm
C}^*(\Delta_\Gamma/\Gamma)$ by bounded linear operators on the
Hilbert space ${\mathcal L}$.
\item The data $({\rm
C}^*(\Delta_\Gamma/\Gamma), \cH, {\mathcal D})$, with the algebra
acting diagonally on ${\mathcal H}={\mathcal L}\oplus
{\mathcal L}$, and the Dirac operator ${\mathcal D}$ of \eqref{D-Lp}
form a spectral triple.
\end{enumerate}
\end{thm}

In \cite{CM1} we recovered arithmetic information such as the local
$L$--factors of \cite{Den2} from the dynamical cohomology and the
data of the spectral triple of Theorem \ref{triple1}.

Recall that, for a curve $X$ over a global field $\mK$, assuming
semi-stability at all places of bad reduction, the local Euler
factor at a place $v$ has the following description
(\cite{Serre}):
\begin{equation}\label{L-factor}
 L_v (H^1(X),s)= \det\left( 1-Fr_v^* N(v)^{-s} | H^1(\bar X,
\Q_\ell)^{I_v} \right)^{-1}.
\end{equation}
Here $Fr_v^*$ is the geometric Frobenius acting on $\ell$-adic
cohomology of $\bar X=X \otimes \Sp(\bar\mK)$, with $\bar\mK$ the
algebraic closure and $\ell$ a prime with $(\ell,q)=1$, where $q$
is the cardinality of the residue field $k(v)$ at $v$. We
denote by $N$ the norm map. The determinant is evaluated on the
inertia invariants $H^1(\bar X, \Q_\ell)^{I_v}$ at $v$ (all of
$H^1(\bar X, \Q_\ell)$ when $v$ is a place of good reduction).

Suppose $v$ is a place of $k(v)$-split degenerate reduction. Then
the completion of $X$ at $v$ is a Mumford curve $X_\Gamma$. In
this case, the Euler factor \eqref{L-factor} takes the following
form:
\begin{equation}\label{loc-factor}
L_v(H^1(X_\Gamma),s)=(1-q^{-s})^{-g}.
\end{equation}
This is computed by the zeta regularized determinant
\begin{equation}\label{det-pm-L}
\det_{\infty,\pi(\cV),i{\mathcal D}}( s ) =
L_v (H^1(X_\Gamma),s)^{-1},
\end{equation}
where
\begin{equation}\label{det-pm}
\det_{\infty,a,i{\mathcal D}}(s) := \exp\left(
-\zeta^\prime_{a,i{\mathcal D},+}(s,0) \right)
\exp\left(-\zeta^\prime_{a,i{\mathcal D},-}(s,0) \right),
\end{equation}
for
\begin{equation}\label{zetapm}
\begin{array}{ll}
\zeta_{a,i{\mathcal D},+} (s,z) := & \sum_{\lambda \in \Sp(i{\mathcal
D})\cap i[0,\infty)}
\Tr(a\Pi_\lambda) (s+\lambda)^{-z} \\[3mm]
\zeta_{a,i{\mathcal D},-} (s,z) := & \sum_{\lambda \in \Sp(i{\mathcal
D})\cap i(-\infty,0)}
\Tr(a\Pi_\lambda) (s+\lambda)^{-z}.
\end{array}\end{equation}
The element $a=\pi({\mathcal V})$ is the projection onto a linear
subspace ${\mathcal V}$ of ${\mathcal H}$, which is obtained
via embeddings of the cohomology of the dual graph
$\Delta_\Gamma/\Gamma$ into the space of cochains of the dynamical
cohomology.

The projection $\pi({\mathcal V})$ acts on the range of the spectral
projections $\hat\Pi_n$ of $D$ as elements $Q_n$ in the AF algebra
core of the ${\rm C}^*$-algebra ${\rm C}^*(\Delta_\Gamma/\Gamma)$.

Notice how, unlike the local factor at infinity, the
factor at the non-archimedean places involves the full spectrum of $D$
and not just its positive or negative part. It is believed that
this difference should correspond to the presence of an underlying
geometric space based on loop geometry, which manifests itself as
loops at the non-archimedean places and as ``half loops'' (holomorphic
disks) at arithmetic infinity.

There is another important difference between the
archimedean and non-archimedean cases. At the archimedean prime  the
local factor is described in terms of zeta functions for a Dirac
operator $D$ (\cf \cite{CM}, \cite{Den}). On the other hand, at the
non-archimedean places, in order to get the correct normalization as
in \cite{Den3}, we need to introduce a rotation of the Dirac operator
by the imaginary unit, $D\mapsto iD$. This rotation corresponds to
the Wick rotation that moves poles on the real line to
poles on the imaginary line (zeroes for the local factor) and appears
to be a manifestation of a rotation from Minkowskian to Euclidean
signature $it \mapsto t$, as already remarked by Manin
(\cite{Man-zeta} p.135), who wrote that {\em ``imaginary time motion''
may be held responsible for the fact that zeroes of $\Gamma(s)^{-1}$ are
purely real whereas the zeroes of all non-archimedean Euler factors
are purely imaginary}. It is expected, therefore, that a more
refined construction would involve a version of spectral
triples for Minkowskian signature.

\vskip .5in

\noindent {\bf Caterina Consani}, Department of Mathematics,
University of Toronto, Canada.

\noindent email: kc\@@math.toronto.edu

\vskip .3in

\noindent {\bf Matilde Marcolli}, Max--Planck--Institut f\"ur
Mathematik, Bonn Germany.

\noindent email: marcolli\@@mpim-bonn.mpg.de

\end{document}